\documentclass[11pt,draft]{article}
\usepackage{amsfonts}
\usepackage{mathrsfs}
\usepackage{amsmath,amsfonts,mathrsfs,amssymb,color}
\usepackage{indentfirst}

\numberwithin{equation}{section}

\newtheorem{theo.}{\quad\, Theorem}[section]
\newtheorem{defi.}{\quad\, Definition}[section]
\newtheorem{lemm.}{\quad\, Lemma}[section]
\newtheorem{coro.}{\quad\, Corollary}[section]

\begin{document}

\title{Existence and uniqueness
of positive solutions for Kirchhoff type beam equations$^*$ }
\author{Jinxiang Wang$^{\ast}$
\\
 {\small  Department of Applied Mathematics, Lanzhou University of Technology, Lanzhou, P.R. Chin}\\
}
\date{} \maketitle
\footnote[0]{E-mail address: wjx19860420@163.com(Jinxiang Wang), \ \ Tel:
86-931-7971297.
%wangh@asu.edu (H. Wang) \
%Telephone: 86-931-7971297
% $^{\ddagger} $xuling$_{-}$216@yahoo.cn
} \footnote[0] {$^*$Corresponding author: Jinxiang Wang.\  Supported by the NSFC(No.11801453). }

 \begin{abstract}
\baselineskip 18pt
This paper is concerned with the existence and uniqueness of positive solution for the fourth order Kirchhoff type problem
 $$
   \left\{\begin{array}{ll}
     u''''(x)-(a+b\int_0^1(u'(x))^2dx)u''(x)=\lambda f(u(x)),\ \ \ \ x\in(0,1),\\
    u(0)=u(1)=u''(0)=u''(1)=0,\\
\end{array} \right.
$$
 where $a>0, b\geq 0$ are constants, $\lambda\in \mathbb{R}$ is a parameter. For the case $f(u)\equiv u$, we use an argument based on the linear eigenvalue problems of fourth order equations and their properties to show that there exists a unique positive solution for all $\lambda>\lambda_{1,a}$, here $\lambda_{1,a}$ is the first eigenvalue of the above problem with $b=0$;  For the case $f$ is sublinear, we prove that there exists a unique positive solution for all $\lambda>0$ and no positive solution for $\lambda<0$ by using bifurcation method.

 \end{abstract}


\vskip 3mm

{\small\bf MR(2000)\ \ \ 34B10, \ 34B18}

\baselineskip 22pt

\section{Introduction}
\vskip 3mm
Consider the following nonlinear fourth order Kirchhoff type problem
$$
   \left\{\begin{array}{ll}
     u''''(x)-(a+b\int_0^1(u'(x))^2dx)u''(x)=\lambda f(u(x)),\ \ \ \ x\in(0,1),\\
    u(0)=u(1)=u''(0)=u''(1)=0,\\
\end{array} \right.\eqno (1.1)
$$
 where $a>0, b\geq0$ are constants, $\lambda\in \mathbb{R}$ is a parameter, $f: \mathbb{R}\rightarrow \mathbb{R}$ is continuous. Due to the presence of the integral term $(b\int_0^1(u'(x))^2dx)u''(x)$, the equation is not a pointwise identity and therefore is a nonlocal integro-differential problem.

 Problem (1.1) describes the bending equilibrium of an extensible beam of length $1$ which is simply supported at two endpoints $x=0$ and $x=1$. The right side term $\lambda f(u)$ in equation stands for a force exerted on the beam by the foundation.
In fact, (1.1) is related to the stationary problem associated with
 $$
 \frac{\partial^{2}u}{\partial t^{2}}+\frac{EI}{\rho A}\frac{\partial^{4}u}{\partial x^{4}}-(\frac{H}{\rho}+\frac{E}{2\rho L}\int^{L}_{0}|\frac{\partial u}{\partial x}|^{2}dx)\frac{\partial^{2}u}{\partial x^{2}}=0,\eqno (1.2)
 $$
which was proposed by Woinowsky-Krieger [1] as a model for the deflection of an extensible beam of length $L$ with hinged ends. In (1.2), $u=u(x,t)$ is the lateral displacement at the space
coordinate $x$ and the time $t$; the letters $H, E, \rho, I$ and $A$ denote, respectively, the tension in the rest
position, the Young elasticity modulus, the density, the cross-sectional moment of inertia and the crosssectional
area. The nonlinear term in the brackets is a correction to the classical Euler-Bernoulli equation
$$
\frac{\partial^{2}u}{\partial t^{2}}+\frac{EI}{\rho A}\frac{\partial^{4}u}{\partial x^{4}}=0,
$$
which does not consider the changes of the tension induced by the variation of the length during the
deflection. This kind of correction was proposed by Kirchhoff [2] to generalize D'Alembert's equation with clamped ends. For this reason (1.1) is often called a Kirchhoff type beam equation. %(Therefore (1) is frequently called a Kirchhoff type equation for a static beam.)
 Other problems involving fourth-order equations of Kirchhoff type can be found in [3-4].

%Equation (1.1) is known as a beam equation of Kirchhoff type and it is based on a model developed in [1].

In the study of problem (1.1) and its generalizations, the nonlocal term under the integral sign causes some mathematical difficulties which make the study of the problem particularly interesting. The existence and multiplicity of solutions for (1.1) and its multi-dimensional case have been studied by several authors, see [5-12] and the references there in. Meanwhile,  numerical methods of (1.1) have been developed in [13-20].

In [5-7], by using variational methods, Ma considered existence and multiplicity of solutions for (1.1) with $\lambda\equiv 1$ under different nonlinear boundary conditions. In [8], based on the fixed point theorems in cones of ordered Banach spaces, Ma studied existence and multiplicity of positive solutions results for (1.1) with right side term $f(x,u,u')$ in equation.

%Based on fixed point theorems in cones, Ma obtained positive solutions of the problem (1.1)( and its generalizations, as well as for equations similar to (1.1))  in [5].  On the topic of (1.1) with other boundary conditions, existence and multiplicity results are also considered in [6]-[8] by using variational methods.
For multi-dimensional case of (1.1) with $\lambda\equiv 1$, Wang et al. studied the existence and multiplicity of nontrivial solutions by using the mountain pass theorem and the truncation method in [9-10]; For a kind of problem similar to (1.1) in $\mathbb{R}^{3}$, Xu and Chen [11] established the existence and multiplicity of negative energy solutions based on the genus properties in critical point theory, and very rencently, Mao and Wang [12] studied the existence of nontrivial mountain-pass type of solutions via the Mountain Pass lemma.

% (Among them, positive solutions were only discussed in [9], for a related problem with m = 0. Our objective is to establish an existence result for positive solutions of (1.1), with m = 0, which is essentially different from the ones in [9]. It is worth noting that in the applications the function m is required to be a line with positive slope. Our main result allows m to have this property. )

%(It should be emphasized that the above mentioned works are of pure theoretical character. In these works the authors only proved the existence of solutions without examples of existing solutions.)
It is worth noticing that, in the above mentioned research work, the uniqueness of solutions for the problem has not been discussed. As far as the author knows,  there are very few results on the uniqueness of solutions for problem (1.1). In [16], when the right side term $\lambda f(u(x))=g(x)$ is nonpositive,  Dang and Luan proved  that problem (1.1) has a unique solution by reducing the problem to a nonlinear equation and proposed an iterative method for finding the solution.  Very recently, by using contraction mapping principle, Dang and Nguyen [17] obtained a uniqueness result for (1.1) in multi-dimensional case with the right side term $\lambda f(u(x))=g(x,u)$ is bounded. To the best of our knowledge, apart from the two works mentioned above, there is no other result on the uniqueness of solutions for nonlocal integro-differential problem (1.1).

Motivated by the above works described, the object of this paper is to study the existence and uniqueness of positive solutions for (1.1),  and our main tool is bifurcation method. It should be emphasized that,  global bifurcation phenomena for fourth order problem (1.1) with $b=0$ have been investigated in [21-23], and [24-28] studied second order Kirchhoff type problem by using the bifurcation theory, but as far as we know, the bifurcation phenomena for fourth order Kirchhoff problem (1.1) has not been discussed.

Concretely, in the present paper we are concerned with problem (1.1) under the two cases: $f(u)\equiv u$ or $f$ is sublinear. For $f(u)\equiv u$, (1.1) can be seen as a nonlinear eigenvalue problem, we use an argument based on the linear eigenvalue problems of fourth order equations and their properties to show that there exists a unique positive solution for all $\lambda>\lambda_{1,a}$, where $\lambda_{1,a}$ is the first eigenvalue of (1.1) with $b=0$; For the case $f$ is sublinear, such as $f(u)=c_{1}u^{p}+c_{2}u^{q}$ ($c_{1},c_{2}\geq 0$, $0<p,q<1$ , see Remark 4.1), we prove that there exists a unique positive solution for all $\lambda>0$ and no positive solution for $\lambda<0$ by using bifurcation method.

The rest paper is arranged as follows:  In Section 2, as preliminaries, we first construct  the operator equation corresponding to (1.1). In Section 3, we deal with the case $f(u)\equiv u$ based on the linear eigenvalue problem of fourth order equations and their properties. Finally, for the case $f$ is sublinear, we discuss the existence and uniqueness of positive solutions for (1.1) by using bifurcation method in Section 4.

\vskip 3mm

\section{Preliminaries}

Let $P:=\{u\in C[0,1]: u(x)\geq 0, \forall\  x\in[0,1]\}$ be the positive cone in $C[0,1]$ and let $U:=P\cup (-P)$.

\noindent{\bf Proposition 2.1}  For each $g\in C[0,1]$, there exists a solution $u$ to the problem $$
   \left\{\begin{array}{ll}
     u''''(x)-(a+b\int_0^1(u'(x))^2dx)u''(x)=g(x),\ \ \ \ x\in(0,1),\\
    u(0)=u(1)=u''(0)=u''(1)=0,\\
\end{array} \right.\eqno (2.1)
$$
and if $g\in U$, then $u$ is unique. Moreover, the operator $T: U\rightarrow U$ defined by
$$T(g):= u$$
is compact.
\noindent

\noindent{\bf Proof.}\  First, when $g\equiv 0$, we prove that (2.1) has only a unique solution $u\equiv 0$. Assume that $u$ is a solution of (2.1)  with $g\equiv 0$, set $w=-u''$, then by (2.1) we have $$
   \left\{\begin{array}{ll}
    -w''(x)+(a+b\int_0^1(u'(x))^2dx)w(x)=0, \ \ \ \ x\in(0,1),\\
    w(0)=w(1)=0,\\
\end{array} \right.\eqno (2.2)
$$
$$
   \left\{\begin{array}{ll}
    -u''(x)=w(x), \ \ \ \ x\in(0,1),\\
    u(0)=u(1)=0.\\
\end{array} \right.\eqno (2.3)
$$We claim that the solution of (2.2) is $w\equiv0$. In fact,  suppose on the contrary that $w\not \equiv 0$ is a solution of (2.2), and without loss of generality, $w(\tau)=\max \{w(x)|x\in[0,1]\}> 0$
for some $\tau\in (0,1)$, then we have $w''(\tau)\leq 0$, which contradicts with $w''(\tau)=(a+b\int_0^1(u'(x))^2dx)w(\tau)>0$.
Substituting $w\equiv0$ in (2.3),  $u\equiv 0$ is an immediate conclusion.

Next, we prove the existence and uniqueness of solutions for (2.1) with $g\neq 0$. For any constant
$R\geq0$, let $u_R$ stands for the unique solution of the linear fourth order problem
 $$
   \left\{\begin{array}{ll}
     u''''(x)-(a+bR)u''(x)=g(x),\ \ \ \ x\in(0,1),\\
    u(0)=u(1)=u''(0)=u''(1)=0,\\
\end{array} \right.\eqno (2.4)
$$
then
$$u_R(x)=\int^{1}_{0}\int^{1}_{0}G_{1}(x,t)G_{2,R}(t,s)g(s)dsdt, \ \ \ \  \ x\in [0,1],\eqno (2.5)
$$
$$u_R''(x)=-\int^{1}_{0}G_{2,R}(x,t)g(t)dt, \ \ \ \  \ x\in [0,1],\eqno (2.6)
$$here $$G_{1}(x,t)=\left\{
\begin{aligned}
t(1-x),\ 0\leq t\leq x\leq1,\\
x(1-t),\ 0\leq x\leq t\leq1,
\end{aligned}
\right.\eqno (2.7)
$$and
  $$G_{2,R}(t,s)=\left\{
\begin{aligned}
\frac{\sinh (\sqrt{a+bR}t) \sinh (\sqrt{a+bR}(1-s))}{\sqrt{a+bR}\sinh \sqrt{a+bR}},\ 0\leq t\leq s\leq1,\\
\frac{\sinh (\sqrt{a+bR}s) \sinh (\sqrt{a+bR}(1-t))}{\sqrt{a+bR}\sinh \sqrt{a+bR}},\ 0\leq s\leq t\leq1,
\end{aligned}
\right.\eqno (2.8)
$$ are Green functions of
$$
   \left\{\begin{array}{ll}
    -u''(x)=0, \ \ \ \ x\in(0,1),\\
    u(0)=u(1)=0,\\
\end{array} \right.\eqno (2.9)
$$
and
$$
   \left\{\begin{array}{ll}
    -w''(t)+(a+bR)w(t)=0, \ \ \ \ t\in(0,1),\\
    w(0)=w(1)=0,\\
\end{array} \right.\eqno (2.10)
$$
respectively. Since $0\leq G_{1}(x,t)\leq G_{1}(x,x)$ and $0 \leq G_{2,R}(t,s)\leq G_{2,R}(t,t)\leq \frac{(\sinh \frac{\sqrt{a}}{2})^{2}}{\sqrt{a}\sinh \sqrt{a}}$, then by (2.5)-(2.8) we have that there exist two positive constants $C_{1}$ and $C_{1}$ such that $$
\|u_{R}\|_{\infty}\leq C_{1}\|g\|_{\infty},\ \ \ \ \ \ \ \ \ \|u_{R}''\|_{\infty}\leq C_{2}\|g\|_{\infty}.\eqno (2.11)
$$
Multiplying the equation in (2.4) by $u_{R}$ and integrating it over $[0,1]$, based on boundary conditions and integration by parts we obtain
$$
\int_{0}^{1} (u_{R}'(x))^{2} dx=\frac{\int_{0}^{1}g(x)u_{R}(x)dx-\int_{0}^{1}(u_{R}''(x))^{2}dx}{a+bR}.\eqno (2.12)
$$
Now to get a solution of (2.1), we only need to find $R$ such that
$$
R=y(R):=\frac{\int_{0}^{1}g(x)u_{R}(x)dx-\int_{0}^{1}(u_{R}''(x))^{2}dx}{a+bR}=\int_{0}^{1} (u_{R}'(x))^{2} dx,\eqno (2.13)
$$that is, find a fixed point of $R=y(R)$. Obviously, $y(0) > 0$. On the other hand, by (2.11) we have$$
|y(R)|=\frac{|\int_{0}^{1}g(x)u_{R}(x)dx-\int_{0}^{1}(u_{R}''(x))^{2}dx|}{a+bR}
\leq \frac{C_{1}\|g\|_{\infty}^{2}+C_{2}^{2}\|g\|_{\infty}^{2}}{a}\leq C.\eqno (2.14)
$$
This concludes the existence of fixed point for $R=y(R)$.

 Now, we show that if $g\in U$, the solution of (2.1) is unique. Without loss of generality, we assume on the contrary that for some $g\in P$,  there exist two solutions $u\neq v$. By (2.5) and (2.6), we have $$u\geq 0,\ \ u''\leq 0; \ \ \ \ \ \ v\geq 0,\ \  v''\leq 0. \eqno (2.15)$$ Since $u$ and $v$ satisfy the equation in (2.1), we have
$$
u''''(x)-v''''(x)-[a+b\int_0^1(u'(x))^2dx](u''(x)-v''(x))
-b[\int_0^1(u'(x))^2dx-\int_0^1(v'(x))^2dx]v''(x)=0.\eqno (2.16)
$$ Set $w=-(u''-v'')$. If  $\int_0^1(u'(x))^2dx=\int_0^1(v'(x))^2dx$, then  (2.2) holds for  $w=-(u''-v'')$ and consequently we can obtain $u\equiv v$ by using similar discussion in Proposition 2.1.
If we assume that $\int_0^1(u'(x))^2dx>\int_0^1(v'(x))^2dx$, then by (2.16) and (2.15) we have
$$
u''''(x)-v''''(x)-[a+b\int_0^1(u'(x))^2dx](u''(x)-v''(x))
\leq0,\eqno (2.17)
$$
that is
$$
-w''(x)+[a+b\int_0^1(u'(x))^2dx]w(x)
\leq0.\eqno (2.18)
$$
We claim that (2.18) implies $w\leq 0$. In fact,  suppose on the contrary that  $w(\tau)=\max \{w(x)|x\in[0,1]\}> 0$
for some $\tau\in (0,1)$, then $w''(\tau)\leq 0$. This contradicts with (2.18) with $x=\tau$. On the other hand, based on boundary conditions and integration by parts, from the assumption $\int_0^1(u'(x))^2dx>\int_0^1(v'(x))^2dx$ we have$$
 \aligned
&\int_0^1(u'(x))^2dx-\int_0^1(v'(x))^2dx\\
&\ \ \ =\int_0^1[u'(x)+v'(x)][u'(x)-v'(x)]dx\\
&\ \ \ =-\int_0^1(u(x)+v(x))(u''(x)-v''(x))dx\\
&\ \ \ =\int_0^1(u(x)+v(x))w(x)dx>0.
\endaligned \eqno (2.19)$$
Since (2.15) guarantees that $u(x)+v(x)\geq 0$, then (2.19) contradicts with $w\leq 0$. The uniqueness of solutions for (2.1) is proved.

At the end, let $T: U\rightarrow C^4[0,1]$ be the operator defined by $Tg=u$, where $u$ is the solution of (2.1). Then by (2.5) and (2.11), we can easily get that $T:U\rightarrow U$ is compact.
\hfill{$\Box$}

\noindent{\bf Remark 2.1 } When $g(x)$ is nonpositive, Dang and Luan [16] proved that problem (2.1) has a unique solution by reducing the problem to a nonlinear equation. Compared with [16], our proof in Proposition 2.1 is more concise.

\vskip 3mm

\section{Nonlinear eigenvalue problem}

In this section, we study (1.1) with $f(u)\equiv u$, that is the nonlinear eigenvalue problem
$$
   \left\{\begin{array}{ll}
     u''''(x)-(a+b\int_0^1(u'(x))^2dx)u''(x)=\lambda u(x),\ \ \ \ x\in(0,1),\\
    u(0)=u(1)=u''(0)=u''(1)=0.\\
\end{array} \right.\eqno (3.1)
$$
The solutions of (3.1) are closely related to the following linear eigenvalue problem:
$$
   \left\{\begin{array}{ll}
     u''''(x)-A u''(x)=\lambda u(x),\ \ \ \ x\in(0,1),\\
    u(0)=u(1)=u''(0)=u''(1)=0.\\
\end{array} \right.\eqno (3.2)
$$
In [29],  Del Pino and Manasevich proposed that: a pair of constants $(\lambda,A)$ such that (3.2) possesses a nontrivial solution will be called an eigenvalue pair,  and the corresponding nontrivial solution will be called an eigenfunction. Furthermore, they proved that the eigenvalue pair $(\lambda,A)$ of (3.2) must satisfy
$$\frac{\lambda}{(k\pi)^{4}}-\frac{A}{(k\pi)^{2}}=1,\ \ \ \text{for some} \ k\in\mathbb{N},$$
and the corresponding eigenfunction is $\varphi_{k}=c \sin k\pi x$($c\neq 0$ is an arbitrary constant).

Now, given a positive constant $A$, we use $\lambda_{1, A}$ to denote the principal eigenvalue of problem (3.2), then we have the following results:

\noindent{\bf Lemma 3.1 }

(i) If $A_{1}, A_{2}$ are positive constants such that $A_{1}< A_{2}$, then $\lambda_{1,A_{1}}< \lambda_{1,A_{2}}$.

(ii) Let $B,C$ be two fixed positive constants. Consider the map $$
\lambda_{1}(\mu):=\lambda_{1,B+\mu C},\ \ \ \ \mu\geq 0,
$$
then $\lambda_{1}(\cdot)$ is a continuous and strictly inreasing function and
$$
\lambda_{1}(0)=\lambda_{1,B},\ \ \ \ \ \ \lim\limits_{\mu\rightarrow +\infty}\lambda_{1}(\mu)=+\infty.
$$

\noindent{\bf Proof. } By [29], we know that the principal eigenvalue $\lambda_{1,A}$  of (3.2) satisfy
$$\frac{\lambda_{1,A}}{\pi^{4}}-\frac{A}{\pi^{2}}=1,\eqno (3.3)$$
and the corresponding first eigenfunction is $\varphi_{1}(x)=c \sin \pi x$, where $c\neq 0$ is an arbitrary constant. According to (3.3),
 $\lambda_{1,A}=(1+\frac{A}{\pi^{2}})\pi^{4},$  then (i) and (ii) are immediate consequences.

 \hfill{$\Box$}

 By using Lemma 3.1, we prove the following results on the nonlinear eigenvalue
problem (3.1):

\noindent{\bf Theorem 3.1 } Problem (3.1) has a positive solution $u_{\lambda}$ if and only if $\lambda\in(\lambda_{1,a},+\infty)$,  moreover, the solution $u_{\lambda}$ is unique and satisfying
$$
\lim\limits_{\lambda\rightarrow \lambda_{1,a}}\|u_{\lambda}\|_{\infty}=0,\ \ \ \ \ \ \ \ \lim\limits_{\lambda\rightarrow +\infty}\|u_{\lambda}\|_{\infty}=+\infty.\eqno (3.4)
$$

\noindent{\bf Proof. }  Assume that $u$ is a positive solution of (3.1), then $\int_0^1(u'(x))^2dx>0$, consequently by Lemma 3.1 (i) we have
$$\lambda=\lambda_{1,a+b\int_0^1(u'(x))^2dx}>\lambda_{1,a}.$$
To any $\lambda\in(\lambda_{1,a},+\infty)$, by Lemma 3.1 (ii),  there exists a unique $t_{0}(\lambda)>0$ such that
$$\lambda_{1,a+bt_{0}}=\lambda,$$ moreover, $$
\lim\limits_{\lambda\rightarrow \lambda_{1,a}}t_{0}(\lambda)=0,\ \ \ \ \ \ \ \lim\limits_{\lambda\rightarrow +\infty}t_{0}(\lambda)=+\infty.\eqno (3.5)
$$
For the fixed $t_{0}$, take appropriate principal eigenfunction $\varphi_{1}(x)=c \sin \pi x(c>0)$ of (3.2) associated to $\lambda_{1,a+bt_{0}}$ such that $$
\int_{0}^{1}(\varphi_{1}'(x))^{2}dx=t_{0}.\eqno (3.6)
$$Then it is easy to see that $u_{\lambda}=\varphi_{1}>0$ is a positive solution of (3.1).

To prove the uniqueness, we assume that there exist two positive solutions $u\neq v$, since
$$\lambda=\lambda_{1,a+b\int_0^1(u'(x))^2dx}
=\lambda_{1,a+b\int_0^1(v'(x))^2dx},$$
then Lemma 3.1 (ii) guarantees that $\int_0^1(u'(x))^2dx=\int_0^1(v'(x))^2dx$ and $u$ is proportional to $v$, which implies that $u=v$.

Finally, we prove (3.4). Since the unique positive solution of (3.1) is  $u_{\lambda}=\varphi_{1}(x)=c_{\lambda} \sin \pi x$, where $c_{\lambda}$ is a positive constant depending on $\lambda$,  then by (3.6) and (3.5), we have $$\lim\limits_{\lambda\rightarrow \lambda_{1,a}}\int_{0}^{1}(u_{\lambda}'(x))^{2}dx
=\lim\limits_{\lambda\rightarrow \lambda_{1,a}}\int_{0}^{1}[(c_{\lambda} \sin \pi x)']^{2}dx=\lim\limits_{\lambda\rightarrow \lambda_{1,a}}\frac{1}{2}c_{\lambda}^{2}\pi^{2} \rightarrow 0,\eqno (3.7)$$and similarly $$\lim\limits_{\lambda\rightarrow +\infty}\int_{0}^{1}(u_{\lambda}'(x))^{2}dx=\lim\limits_{\lambda\rightarrow +\infty}\frac{1}{2}c_{\lambda}^{2}\pi^{2}\rightarrow \infty,\eqno (3.8) $$
that is
$$
\lim\limits_{\lambda\rightarrow \lambda_{1,a}}c_{\lambda} \rightarrow 0,\ \ \ \ \ \ \ \lim\limits_{\lambda\rightarrow +\infty}c_{\lambda}\rightarrow +\infty, \eqno (3.9)
$$
then (3.4) is an immediate consequence.\hfill{$\Box$}

\vskip 3mm

\section{The sublinear case}

In this section, we study (1.1) when the nonlinear term $f$ is sublinear which means that $f$ satisfying:

\noindent(H1) $f: \mathbb{R}\rightarrow \mathbb{R}$ is continuous, $f(s)>0$ for all $s>0$, $f(0)=0$ and $
f_{0}:=\lim\limits_{s\rightarrow 0+}\frac{f(s)}{s}=+\infty;
$

\noindent(H2) $
f_{\infty}:=\lim\limits_{s\rightarrow +\infty}\frac{f(s)}{s}=0.
$

 The main tool we will use in this section is global bifurcation theory. %theorem for mappings which are not necessary smooth.

%\noindent{\bf Lemma 4.1} \ (Rabinowitz [30]). Let $X$ be a real reflexive Banach space. Let $T :  \mathbb{R}\times X \mapsto X$ be completely continuous such that $T(\lambda,0)=0,\ \forall \ \lambda\in \mathbb{R}$. Let $l_{1}, l_{2}\in \mathbb{R} (l_{1}<l_{2})$ be such that $u=0$ is an isolated solution of the equation$$
%u-T(\lambda,u)=0,\ u\in X \eqno (4.1)
%$$for $\lambda=l_{1}$ and $\lambda=l_{2}$, where $(l_{1},0), (l_{2},0)$ are not bifurcation points of (4.1). Furthermore, assume that
%$$
%deg(I-T(l_{1}, \cdot),B_{\rho}(0),0)\neq deg(I-T(l_{2}, \cdot),B_{\rho}(0),0),
%$$
%where $B_{\rho}(0)$ is an isolating neighborhood of the trivial solution.
%Let
%$$
%\mathcal{S}:=\overline{\{(\lambda,u): (\lambda,u)\  \text{is a solution of (4.1) with}\ u\neq 0 \}}\cup ([l_{1},l_{2}]\times \{0\}).
%$$
%Then there exists a connected component $\mathcal{C}$ of $\mathcal{S}$ containing $[l_{1},l_{2}]\times \{0\}$, and either

%(i) $\mathcal{C}$ is unbounded in $\mathbb{R}\times X$, or

%(ii) $\mathcal{C} \cap [(\mathbb{R}\setminus[l_{1},l_{2}])\times \{0\}]\neq \emptyset$

%To apply Lemma 4.1,

 We first state some notations. Let $X:=\{u\in C^{2}[0,1]:u(0)=u(1)=u''(0)=u''(1)=0\}$ with the norm $\|u\|_{X}=\max\{\|u\|_{\infty},\|u'\|_{\infty},\|u''\|_{\infty}\}$. $B_{\rho}:=\{u\in X:\|u\|_{X}< \rho\}$. For any $u\in X$,  denote $u^{+}=\max\{u,0\}$. Define the operator $F: \mathbb{R}\times X \mapsto X$ by
 $$
F(\lambda,u)(x):=T(\lambda f(u^{+}(x))), \eqno (4.1)
$$
where $T$ is the operator defined in Proposition 2.1. Then it is easy to see that $u$ is a nonnegative solution of (1.1) if and only if
$$
u=F(\lambda,u). \eqno (4.2)
$$
 Since the map from $X$ into $U:= P \cup(-P)$ defined by $u\mapsto \lambda f(u^{+})$ is continuous, and $C^{4}[0,1]\cap X$ is compactly imbedded in $X$, then by Proposition 2.1, the operator $F :  \mathbb{R}\times X \mapsto X$ as in (4.1) is completely continuous. In order to prove the main result of this section, we need the following lemmas.

\noindent{\bf Lemma 4.1 } For any fixed $\lambda< 0$, there exists a number $\rho>0$ such that $$deg(I-F(\lambda, \cdot),B_{\rho}(0),0)=1.$$

\noindent{\bf Proof. } First, we claim that there exists $\delta>0$ such that
$$
u\neq tF(\lambda,u)=tT(\lambda f(u^{+})) \ \ \ \text{for all}\ \ u\in \overline{B_{\delta}},\  u\neq0 \ \ \text{and} \ t\in [0,1].
$$
Suppose on the contrary that there exist sequence $\{u_{n}\}$ in  $X\setminus 0$ with $\|u_{n}\|_{X}\longrightarrow 0$ and $\{t_{n}\}$ in $[0,1]$ such that
$$u_{n}=t_{n}F(\lambda,u_{n})=t_{n}T(\lambda f(u_{n}^{+})),$$
that is
$$u_{n}''''(x)-(a+b\int_0^1(u_{n}'(x))^2dx)u_{n}''(x)=t_{n}\lambda f(u_{n}^{+}(x))\leq 0,\ \ \ \ x\in(0,1).\eqno (4.3)
$$
Set $w_{n}=-u''_{n}$, then by (4.3) we can get an inequality for $w_{n}$ similar to (2.18), which can deduce that $w_{n}\leq0$. Consequently,  $-u''_{n}=w_{n}\leq 0$ and $u_{n}(0)=u_{n}(1)=0$ guarantee that $u_{n}\leq 0$, which implies $f(u_{n}^{+})\equiv 0$ according to (H1). Then by Proposition 2.1, (4.3) has only a unique solution $u_{n}\equiv 0$, a contradiction with $u_{n}\in X\setminus 0$.

Take $\rho\in (0,\delta]$, according to the homotopy invariance of topological degree, we have
$$
deg(I-F(\lambda, \cdot),B_{\rho}(0),0)=deg(I,B_{\rho}(0),0)=1.
$$\hfill{$\Box$}

\noindent{\bf Lemma 4.2} \ For any fixed $\lambda> 0$, there exists a number $\rho>0$ such that $$deg(I-F(\lambda, \cdot),B_{\rho}(0),0)=0.$$

\noindent{\bf Proof. } First, take a $\psi \in X, \psi>0,$ we claim that there exists $\delta>0$ such that
$$
u\neq T(\lambda f(u^{+})+t\psi) \ \ \ \text{for all}\ \ u\in \overline{B_{\delta}},\  u\neq0 \ \ \text{and} \ t\in [0,1].
$$
Suppose on the contrary that there exist sequence $\{u_{n}\}$ in  $X\setminus 0$ with $\|u_{n}\|_{X}\longrightarrow 0$ and $\{t_{n}\}$ in $[0,1]$ such that
$$u_{n}= T(\lambda f(u_{n}^{+})+t_{n}\psi),$$
that is
$$u_{n}''''(x)-(a+b\int_0^1(u_{n}'(x))^2dx)u_{n}''(x)=\lambda f(u_{n}^{+})+t_{n}\psi ,\ \ \ \ x\in(0,1),\eqno (4.4)
$$
Since $t_{n}\psi> 0$, from the similar argument in Lemma 4.2 we have that $u_{ n} >0$.

%On the other hand, multiplying (4.4) by $u_{n}$ and integrating it over $[0,1]$, based on boundary conditions and integration by parts we obtain
%$$
%\int_{0}^{1} (u_{n}'(x))^{2} dx=\frac{\int_{0}^{1}(\lambda f(u_{n}^{+})+t_{n}\psi(x))u_{n}(x)dx-\int_{0}^{1}(u_{n}''(x))^{2}dx}
%{a+b\int_0^1(u_{n}'(x))^2dx}.
%$$
On the other hand, $\|u_{n}\|_{X}\longrightarrow 0$ implies that %, then $\lambda f(u_{n}^{+})+t_{n}\psi(x)$ is bounded, and by using the same argument in (2.11) and (2.14) we can get
$$\int_{0}^{1} (u_{n}'(x))^{2} dx\leq C$$
for some positive constant $C$. Hence, according to Lemma 3.1 we have that
$$
\lambda_{1,a+b\int_0^1(u_{n}'(x))^2dx}\leq \lambda_{1,a+b C}:=\Lambda.
$$
Fix this value of $\Lambda$, since $\|u_{n}\|_{\infty}\longrightarrow 0$, then according to (H1), for $n$ large we have that $\lambda f(u_{n}^{+})>\Lambda u_{n}$. Combining this with $u_{n}''\leq 0$ we can get
$$u_{n}''''(x)-(a+bC)u_{n}''(x)\geq u_{n}''''(x)-(a+b\int_0^1(u_{n}'(x))^2dx)u_{n}''(x) =\lambda f(u_{n}^{+})+t_{n}\psi>\Lambda u_{n},
$$
which implies that $\lambda_{1,a+bC}>\Lambda$, a contradiction!

Take $\rho\in (0,\delta]$, according to the homotopy invariance of topological degree, we have
$$
deg(I-F(\lambda, \cdot),B_{\rho}(0),0)=deg(I-T(\lambda f(\cdot)+\psi),B_{\rho}(0),0)=0.
$$\hfill{$\Box$}

Now, we are ready to consider the bifurcation of positive solutions of (1.1) from the line of trivial solutions $\{(\lambda,0)\in \mathbb{R}\times X: \lambda\in\mathbb{R}\}$.
%In the following result, we show that from the trivial solution emanates an unbounded continuum of positive solution.

\noindent{\bf Theorem 4.1 } Assume that (H1) and (H2) hold. Then from $(0,0)$ there emanate an unbounded continuum $\mathcal{C}_{0}$ of positive solutions of (1.1) in $\mathbb{R}\times X$.

% we will apply the following global bifurcation theorems for mappings which are not necessary smooth to get a global description of the branches of positive solutions of (1.1).

\noindent{\bf Proof of Theorem 4.1.} By an argument similar to that of [30, Proposition 3.5], using Lemma 4.1 and 4.2, we can show that $(0,0)$ is a bifurcation point from the line of trivial solutions $\{(\lambda,0)\in \mathbb{R}\times X: \lambda\in\mathbb{R}\}$ for the equation (4.2), and there exists a connected component $\mathcal{C}_{0}$ of positive solutions of (4.2) containing $(0,0)$, either

(i) $\mathcal{C}_{0}$ is unbounded in $\mathbb{R}\times X$, or

(ii) $\mathcal{C}_{0} \cap [\mathbb{R}\setminus0\times \{0\}]\neq \emptyset$.

To prove the unboundedness of $\mathcal{C}_{0}$ , we only need to show that the case (ii)  cannot occur, that is: $\mathcal{C}_{0}$ can not meet $(\lambda,0)$ for any $\lambda\neq 0$.   It is easy to see that for $\lambda<0$ problem (1.1) does not possess a positive solution. For the case $\lambda>0$, we assume on the contrary that there exist some $\lambda_{0}>0$ and a sequence of parameters $\{\lambda_{n}\}$ and corresponding positive solutions $\{u_{n}\}$ of (1.1) such
that $\lambda_{n}\longrightarrow \lambda_{0}$ and $\|u_{n}\|_{X} \longrightarrow 0$. Since $\|u_{n}\|_{\infty} \longrightarrow 0$, then by (H1), for fixed $\varepsilon\in(0,\lambda_{0})$ there exists $n_{0}\in \mathbb{N}$ such that when $n> n_{0}$ we have
$$u_{n}''''(x)-(a+b\int_0^1(u_{n}'(x))^2dx)u_{n}''(x)=\lambda_{n} f(u_{n})\geq (\lambda_{0}-\varepsilon)f(u_{n})>\Lambda u_{n},
$$
where $\Lambda$  is defined as in Lemma 4.2. Now, we can get a contradiction in a similar way that in the proof of Lemma 4.2. \hfill{$\Box$}

The main result of this section is following:

\noindent{\bf Theorem 4.2} \  Assume that (H1) and (H2) hold, then (1.1) has a positive solution if and only if $\lambda>0$. In addition, if $f$ is monotone increasing and there exists $\alpha\in (0,1)$ such that $$
f(\tau s)\geq \tau^{\alpha}f(s)\eqno (4.5)
$$ for any $\tau\in(0,1)$ and $s>0$, then the positive solution of (1.1) is unique.

\noindent{\bf Proof.} By Theorem 4.1, there exists an unbounded continuum $\mathcal{C}_{0}\in \mathbb{R}\times X$ of positive solutions of (1.1).
We will show that $\|u\|_{X}$ is bounded for any fixed $\lambda>0$, that is, $\mathcal{C}_{0}$ can not blow up at finite $\lambda\in (0,+\infty)$.
 %that is $$\sup\{\lambda|\  (\lambda,u)\in \mathcal{C}_{0}\}=\infty$$
 To do this, we first prove $\|u\|_{\infty}$ is bounded for any fixed $\lambda>0$. Assume on the contrary that there exist $\lambda_{0}>0$ and a sequence of parameters $\{\lambda_{n}\}$ and corresponding positive solutions $\{u_{n}\}$ of (1.1) such
that $\lambda_{n}\longrightarrow \lambda_{0}$, $\|u_{n}\|_{\infty}\longrightarrow \infty$. %We first claim that $\|u_{n}\|_{X}\longrightarrow \infty$ implies $\|u_{n}\|_{\infty}\longrightarrow \infty$. In fact,
Since
$$u_{n}''''(x)-(a+b\int_0^1(u_{n}'(x))^2dx)u_{n}''(x)=\lambda_{n} f(u_{n}),\eqno (4.6)
$$
divide (4.6) by $\|u_{n}\|_{\infty}$ and set $v_{n}=\frac{u_{n}}{\|u_{n}\|_{\infty}}$, then we get
$$
v_{n}''''(x)-(a+b\int_0^1(u_{n}'(x))^2dx)v_{n}''(x)=\lambda_{n} \frac{f(u_{n}(x))}{\|u_{n}\|_{\infty}}.\eqno (4.7)
$$
Multiplying (4.7) by $v_{n}$ and integrating it over $[0,1]$, based on boundary conditions and integration by parts we obtain
$$
\int_{0}^{1} (v_{n}'(x))^{2} dx=
\frac{\int_{0}^{1}\lambda_{n}\frac{f(u_{n}(x))}{\|u_{n}\|_{\infty}}v_{n}(x)dx-\int_{0}^{1}(v_{n}''(x))^{2}dx}{a+b\int_0^1(u_{n}'(x))^2dx}.\eqno (4.8)
$$
Since $\|v_{n}\|_{\infty}\equiv 1$, $\{\lambda_{n}\}$ is bounded and (H2) guarantees that $
\frac{f(u_{n}(x))}{\|u_{n}\|_{\infty}}\longrightarrow 0$ as $n\longrightarrow \infty$, then  (4.8) implies
$$
0\leq\int_{0}^{1} (v_{n}'(x))^{2} dx\leq
\frac{\int_{0}^{1}\lambda_{n}\frac{f(u_{n}(x))}{\|u_{n}\|_{\infty}}v_{n}(x)dx
}{a}\longrightarrow 0\ \ \ \ \text{as}\ n\longrightarrow \infty,
$$
that is $\|v_{n}'\|_{\infty}\longrightarrow 0$.  By the boundary conditions $v_{n}(0)=v_{n}(1)=0$, there exist $\xi_{n}\in (0,1)$ such that $v_{n}(x)=\int_{\xi_{n}}^{x}v_{n}'(t)dt,\ \forall x\in [0,1]$. Combining this with $\|v_{n}'\|_{\infty}\longrightarrow 0$ we can conclude that $\|v_{n}\|_{\infty}\longrightarrow 0$. This contracts with $\|v_{n}\|_{\infty}\equiv1$, and then we get the boundedness of $\|u\|_{\infty}$. Next, we show that the boundedness of $\|u\|_{\infty}$ can deduce the boundedness of $\|u'\|_{\infty}$ and $\|u''\|_{\infty}$.
Since
$$u''''(x)-(a+b\int_0^1(u'(x))^2dx)u''(x)=\lambda f(u(x)),\eqno (4.9)
$$
multiplying (4.9) by $u$ and integrating it over $[0,1]$, similarly we can obtain
$$
\int_{0}^{1} (u'(x))^{2} dx=
\frac{\int_{0}^{1}\lambda f(u(x))u(x)dx-\int_{0}^{1}(u''(x))^{2}dx}{a+b\int_0^1(u'(x))^2dx}\leq
\frac{\int_{0}^{1}\lambda f(u(x))u(x)dx
}{a}.\eqno (4.10)
$$
(4.10) implies that $\|u'\|_{\infty}$ is bounded, and consequently, $\|u''\|_{\infty}$ is bounded too. According to the definition of $\|u\|_{X}$, the above conclusion means that $\|u\|_{X}$ is bounded for any fixed $\lambda>0$. Combining this with the unboundedness of $\mathcal{C}_{0}$, we conclude that $\sup\{\lambda|\ (\lambda,u)\in\mathcal{C}_{0}\}=\infty$, then for any $\lambda> 0$ there exists a positive solution for (1.1).

Now, we prove that if $f$ is monotone increasing and satisfies (4.5), then (1.1) has only a unique positive solution. Assume that there exist two positive solutions $u\neq v$ corresponding to some fixed $\lambda>0$. If  $\int_0^1(u'(x))^2dx=\int_0^1(v'(x))^2dx=R>0$, consider the problem
$$
   \left\{\begin{array}{ll}
     \omega''''(x)-(a+b\int_0^1(u'(x))^2dx)\omega''(x)=\lambda f(\omega(x)),\ \ \ \ x\in(0,1),\\
    \omega(0)=\omega(1)=\omega''(0)=\omega''(1)=0,\\
\end{array} \right.\eqno (4.11)
$$
and its corresponding integral operator $H: P\rightarrow P $  given by
$$
H(\omega)=T(\lambda f(\omega))=\lambda\int^{1}_{0}\int^{1}_{0}G_{1}(x,t)G_{2,R}(t,s)f(\omega(s))dsdt.
$$
By the monotonicity of $f$ and (4.5), $H$ is an increasing $\alpha$-concave operator according to [31, Definition 2.3], then by [31, Theorem 2.1, Remark 2.1], the operator equation $H(\omega)=\omega$  has a unique solution, which is also the unique positive solution of (4.11). That is, $u=v.$

If we assume that $\int_0^1(u'(x))^2dx>\int_0^1(v'(x))^2dx$, since $v''\leq0$, we have
$$
v''''(x)-(a+b\int_0^1(u'(x))^2dx)v''(x)\geq  v''''(x)-(a+b\int_0^1(v'(x))^2dx)v''(x)=\lambda f(v(x)), \eqno (4.12)
$$
which means that $v$ is actually an upper solution of (4.11).  Constructing an iterative sequence $v_{n+1}=Hv_{n},n=0,1,2,\ldots$, where $v_{0}=v$, then (4.12) and the monotonicity of $f$ guarantee that $\{v_{n}\}$ is decreasing.  Moreover, by [31, Theorem 2.1, Remark 2.1], $\{v_{n}\}$ must converge to the unique solution  $u$ of (4.11), %according to [31, Definition 4]. Obviously, problem (4.9) has a lower solution $0$. Then, according to [31,Theorem 7], there exists a solution $\omega$ for problem (4.9) satisfying
%$$
%0\leq\omega(x)\leq v(x), \ \ \ \ \forall x\in [0,1].
%$$
%Since the solution of problem (4.9) is unique,  then $\omega$ and $u$ are actually the same, that is, $\omega \equiv u$,
and consequently we have
$$
0\leq u(x)\leq v(x), \ \ \ \ \forall x\in [0,1].\eqno (4.13)
$$
On the other hand, based on boundary conditions and integration by parts, from the assumption $\int_0^1(u'(x))^2dx>\int_0^1(v'(x))^2dx$ we have that $$
 \aligned
&\int_0^1(u'(x))^2dx-\int_0^1(v'(x))^2dx\\
&\ \ \ =\int_0^1[u'(x)+v'(x)][u'(x)-v'(x)]dx\\
&\ \ \ =-\int_0^1(u(x)-v(x))(u''(x)+v''(x))dx>0,
\endaligned \eqno (4.14)$$
since $-(u''(x)+v''(x))\geq 0$ following from (2.15), then (4.14) contradicts with (4.13). This concludes the proof. \hfill{$\Box$}

\noindent{\bf Remark 4.1 }  Let $c_{1},c_{2}$ are nonnegative constants satisfying $c_{1}^{2}+c_{2}^{2}\neq 0$, $0<p,q<1$, then it is easy to check that the function $$
f(u)=c_{1}u^{p}+c_{2}u^{q}
$$
is increasing and satisfies (H1),(H2) and (4.5).  Consequently,  Theorem 4.2 guarantees that the problem $$
   \left\{\begin{array}{ll}
     u''''(x)-(a+b\int_0^1(u'(x))^2dx)u''(x)=\lambda (c_{1}u^{p}+c_{2}u^{q}),\ \ \ \ x\in(0,1),\\
    u(0)=u(1)=u''(0)=u''(1)=0,\\
\end{array} \right.
$$
has a positive solution if and only if $\lambda>0$, moreover,  the positive solution is unique.

\vskip 3mm

\vskip 5mm

\noindent{\bf Acknowledgements}

\noindent The authors are very grateful to the anonymous referees for their valuable
suggestions. This work was supported by the NSFC (No. 11801453).

\vskip 12mm

\centerline {\bf REFERENCES}\vskip5mm\baselineskip 0.45cm
\begin{description}
\baselineskip 15pt

\item{[1]}~  S. Woinowsky-Krieger, The effect of axial force on the vibration of hinged bars, J. Appl. Mech. 17 (1950) 35-36.

\item{[2]}  G. Kirchhoff, Vorlesungen ¨¹ber Mathematiche Physik: Mechanik, Teubner, Leipzig, 1876.

\item{[3]} K. Doppel, K. Pfl$\ddot{u}$ger and W. Herfort, A nonlinear beam equation arising in the theory of bluff bodies, Z. Anal. Anwend. 16, 945-960,(1997).

\item{[4]} L. A. Medeiros, On a new class of nonlinear wave equations, J. Math. Anal. Appl. 69, 252-262, (1979).

\item{[5]}  T. F. Ma, Existence results for a model of nonlinear beam on elastic bearings, Appl. Math. Lett. 13 (5) (2000) 11-15.

\item{[6]}  T. F. Ma, Existence results and numerical solutions for a beam equation with nonlinear boundary conditions, Appl. Numer.
Math. 47 (2003) 189-196.

\item{[7]}  T. F. Ma, Positive solutions for a nonlinear Kirchhoff type beam equation, Appl. Math. Lett. 18 (2005) 479-482.

\item{[8]} ~T. F. Ma, Positive solutions for a nonlocal fourth order equation of Kirchhoff type,  Discrete Contin Dyn Syst, 2007, (Supplement): 694-703.

\item{[9]} F. L. Wang, Y. K. An, Existence and multiplicity of solutions for a fourth-order elliptic equation, Bound. Value Probl. 2012 (2012) 6.

\item{[10]}  F. L. Wang, M. Avci, Y. K. An, Existence of solutions for fourth order elliptic equations of Kirchhoff type, J. Math. Anal. Appl. 409 (2014) 140-146.

\item{[11]} L. P. Xu, H. B. Chen, Multiplicity results for fourth order elliptic equations of Kirchhoff-type, Acta Mathematica Scientia, Volume 35, Issue 5,  2015,  1067-1076.

\item{[12]}  A. M. Mao, W. Q. Wang, Nontrivial solutions of nonlocal fourth order elliptic equation of Kirchhoff type in $\mathbb{R}^{3}$, J. Math. Anal. Appl. 459 (2018) 556-563.

%\item{[12]} L. P. Xu, H. B. Chen, Existence and multiplicity of solutions for fourth-order elliptic equations of Kirchhoff type via genus theory. Boundary Value Problems 2014, 2014: 212

%\item{[13]} M. Ferrara, S. Khademloo, S. Heidarkhani, Multiplicity results for perturbed fourth-order Kirchhoff type elliptic problems,  Applied Mathematics and Computation 234 (2014) 316-325.

\item{[13]} J. Y. Shin, Finite-element approximations of a fourth-order differential equation, Comput. Math. Appl. 35 (8) (1998) 95-100.

\item{[14]} M. R. Ohm , H. Y. Lee , J. Y. Shin , Error estimates of finite-element approximations for a fourth-order differential equation, Comput. Math. Appl. 52 (2006) 283-288.

\item{[15]} J. Peradze , A numerical algorithm for a Kirchhoff-type nonlinear static beam, J. Appl. Math. 12 (2009) ID 818269.

\item{[16]} Q. A. Dang, V. T. Luan, Iterative method for solving a nonlinear fourth order boundary value problem, Comput. Math. Appl. 60 (2010) 112-121.

\item{[17]} Q. A. Dang, T. H. Nguyen, Existence results and iterative method for solving a nonlinear biharmonic equation of Kirchhoff type,  Comput. Math. Appl. 76 (2018) 11-22.

\item{[18]} H. Temimi , A. R. Ansari , A. M. Siddiqui , An approximate solution for the static beam problem and nonlinear integro-differential equations, Comput. Math. Appl. 62 (2011) 3132-3139.

\item{[19]} Q. Q. Zhuang, Q. W. Ren, Numerical approximation of a nonlinear fourth-order integro-differential equation by spectral method, Appl. Math. Comput. 232 (2014) 775-783

\item{[20]}  Q. W. Ren, H. J. Tian, Numerical solution of the static beam problem by Bernoulli collocation method, Appl. Math. Model. 40 (2016) 8886-8897.

\item{[21]}~ R. Y. Ma, Existence of positive solutions of a fourth-order boundary value problem, Appl. Math. Comput. 168 (2005), 1219-1231.

\item{[22]}~ P. Korman, Uniqueness and exact multiplicity of solutions for a class of fourth-order semilinear problems, Proc. Roy. Soc. Edinburgh Sect. A 134(1) (2004), 179-190.

\item{[23]}~ B. P. Rynne, Global bifurcation for $2m$th-order boundary
value problems and infinitely many solutions of superlinear
problems, J. Differential Equations 188 (2003), 461-472.

\item{[24]} Z. P. Liang, F. Y. Li, J. P. Shi, Positive solutions of Kirchhoff-type non-local elliptic equation: a bifurcation
approach, Proc. Roy. Soc. Edinburgh Sect. A 147 (2017), 875-894.

\item{[25]} G. M. Figueiredo, et al., Study of a nonlinear Kirchhoff equation with non-homogeneous material, J. Math. Anal. Appl. 416 (2014) 597-608.

\item{[26]} Z. P. Liang, F. Y. Li, J. P. Shi, Positive solutions to Kirchhoff type equations with nonlinearity having prescribed asymptotic behavior, Ann. Inst. H. Poincar¨¦
Anal. Non Lin¨¦aire 31 (2014) 155-167.

\item{[27]} G. W. Dai, H. Y. Wang, B. X. Yang, Global bifurcation and positive solution for a class of fully
nonlinear problems,  Comput. Math. Appl. 69 (2015) 771-776.

\item{[28]} A. Ambrosetti, D. Arcoya, Positive solutions of elliptic Kirchhoff equations, Adv. Nonlinear Stud. 17(1), (2017) 3-15.

\item{[29]} M. A. Del Pino, R. F. Manasevich, Existence for a fourth-order boundary value problem under a two-parameter
nonresonance condition, Proc. Amer. Math. Soc. 112 (1991) 81-86.

%\item{[30]} P. H. Rabinowitz, Some aspects of nonlinear eigenvalue problems, in: Rocky Mountain Consortium Symposium on Nonlinear Eigenvalue Problems (N.M. Santa Fe, 1971). Rocky Mountain J. Math. 3 (1973) 161-202.

\item{[30]}  A. Ambrosetti, P. Hess, Positive solutions of asymptotically linear elliptic eigenvalue problems, J. Math. Anal. Appl. 73 (1980) 411-422.

\item{[31]}  C. B. Zhai, C. R. Jiang, S. Y. Li, Approximating monotone positive solutions of a nonlinear fourth-order boundary value problem via sum operator method, Mediterr. J. Math.  14 (2017) 77.

%\item{[32]} X. N. Lin, D. Q. Jiang, X. Y. Li, Existence and uniqueness of solutions for singular fourth-order boundary value problems, J.  Comput. Appl. Math. 196 (2006) 155-161.

%\item{[32]} R. Vrabel, On the lower and upper solutions method for the problem of elastic beam with hinged ends, J. Math. Anal. Appl. 421(2), 1455-1468 (2015)

\end{description}
\end{document}